\documentclass[12pt]{amsart}
\usepackage{amssymb,latexsym,comment,url}

\newcommand{\CC}{\operatorname{\mathcal C}}
\newcommand{\PP}{{\mathbb P}}
\newcommand{\Q}{{\mathbb Q}}

\newcommand{\Z}{{\mathbb Z}}

\newenvironment{Proof}{\par\noindent{\sc Proof:}}%
                      {\hspace*{\fill}\nobreak$\Box$\par\medskip}
                       {\hspace*{\fill}\nobreak$\Box$\par\medskip}

\newtheorem{Proposition}{Proposition}[section]
\newtheorem{Theorem}[Proposition]{Theorem}
\newtheorem{Lemma}[Proposition]{Lemma}
\newtheorem{Corollary}[Proposition]{Corollary}

\theoremstyle{definition}

\newtheorem{Remark}[Proposition]{Remark}
 
\renewcommand{\baselinestretch}{1.1}

\begin{document}

\title[Large $F$-Diophantine sets]%
{On large $F$-Diophantine sets}

\author[M. Sadek]%
{Mohammad~Sadek}
\address{American University in Cairo, Mathematics and Actuarial Science Department, AUC Avenue, New Cairo, Egypt}
\email{mmsadek@aucegypt.edu}

\author[N. El-Sissi]%
{Nermine El-Sissi}
\address{American University in Cairo, Mathematics and Actuarial Science Department, AUC Avenue, New Cairo, Egypt}
\email{nelsissi@aucegypt.edu}

\let\thefootnote\relax
\begin{abstract}
Let $F\in\Z[x,y]$ and $m\ge2$ be an integer. A set $A\subset \Z$ is called an $(F,m)$-Diophantine set if $F(a,b)$ is a perfect $m$-power for any $a,b\in A$ where $a\ne b$. If $F$ is a bivariate polynomial for which there exist infinite $(F,m)$-Diophantine sets, then there is a complete qualitative characterization of all such polynomials $F$. Otherwise, various finiteness results are known. We prove that given a finite set of distinct integers $ S$ of size $n$, there are infinitely many bivariate polynomials $F$ such that $ S$ is an $(F,2)$-Diophantine set. In addition, we show that the degree of $F$ can be as small as $\displaystyle 4\lfloor n/3\rfloor$.
\end{abstract}

\maketitle


\section{Introduction}
The polynomial $F(x,y)=xy+1$ is a bivariate polynomial of degree $2$ that satisfies the following property: there are infinitely many sets of integers of size $4$, $\displaystyle \left\{a_1,a_2,a_3,a_4\right\}$, such that $F(a_i,a_j)$ is a perfect square for any $1\le i,j\le 4$ where $i\ne j$. While it is known that there are no such sets of size $6$, it was proved that there are at most finitely many such sets of size $5$, see \cite{Dujella1}. Recently, a proof of the non-existence of such sets of integers of size $5$ was presented, \cite{He}. In fact, a set of integers $S$ of size $m$ enjoying the property that $F(a_i,a_j)$, $a_i,a_j\in S$, $i\ne j$, is a perfect square is said to be a Diophantine $m$-tuple. So the previous facts can be restated as follows: there are infinitely many Diophantine quadruples while there are no Diophantine $m$-tuples when $m\ge 5$.

One may replace the polynomial $F(x,y)=xy+1$ with any bivariate polynomial in $\Z[x,y]$ and then look for sets of distinct integers $S$ such that $F(a_i,a_j)$ is a perfect $m$-power, $m\ge 2$, for any $a_i,a_j\in S$, $i\ne j$. The set $S$ is then called an $(F,m)$-Diophantine set. One remarks that this definition allows patterns that were not seen when tackling the Diophantine $m$-tuple problem. For example, there is an infinite $(F,m)$-Diophantine set when $F(x,y)=xy$ for any $m\ge 2$.

Many questions may be posed now. Fixing a bivariate polynomial $F(x,y)\in\Z[x,y]$, is it possible that there is an infinite $(F,m)$-Diophantine set? If not, then what is the size of the largest $(F,m)$-Diophantine set? Or at least, can one find an upper bound on the size of such set? In addition, given a set of integers $S$, can we construct a bivariate polynomial $F$ such that $S$ is an $(F,m)$-Diophantine set? If so, then what is the smallest possible degree of $F$?

$(F,m)$-Diophantine sets were introduced in \cite{BercczesDujella}. The polynomials $F(x,y)$ for which there are infinite $(F,m)$-Diophantine sets $S$ are completely characterized. They are of the form $f(x)^{s}g(y)^t$, for some $s,t\ge0$, where $0\le \deg f,\deg g\le 2$, see \cite[Theorem 2.2]{BercczesDujella}. Otherwise, bounds on the size of $S$ are given. Those bounds are not explicit and they depend on $F$ and $m$. Finally, the authors gave explicit bounds on the size of $(F,2)$-Diophantine sets when $F$ is a quartic bivariate polynomial of certain form.

In this note, given a set of distinct integers $\displaystyle S=\{x_1,\ldots, x_n\}$, we prove that there are infinitely many polynomials $f(x)\in\Z[x]$ such that $S$ is an $(F_f,2)$-Diophantine set, where $F_f(x,y)=f(x)f(y)$. Moreover, the degree of $f(x)$ is $n-2$. In order to find such $f$, we use ideas of Yamagishi, \cite{Yamagishi1,Yamagishi2,Yamagishi3}. We associate a hypersurface in the projective space to the set $S$. We then prove that the latter variety is rational, hence it possesses infinitely many rational points. The rationality of the latter variety implies the existence of the polynomial $f$. In particular, allowing $F$ to run all over bivariate polynomials with integer coefficients, we prove that $(F,2)$-Diophantine sets can be arbitrarily large.

A further refinement on the degree of $f$ is displayed. We associate a determinantal variety consisting of the intersection of diagonal quadrics to $S$. By observing that this variety has several rational points, we prove that it is rational. This implies the existence of infinitely many polynomials of degree $\displaystyle 2 \lfloor  \#S/3\rfloor$ such that $S$ is an $(F_f,2)$-Diophantine set, where $\lfloor\textrm{ . }\rfloor$ is the greatest integer function.

\section{The polynomial $F(x,y)$}

Given a polynomial $f(x)\in\Z[x]$, we will write $F_f(x,y)$ for $f(x)f(y)$. In this section, we will show that for any positive integer $k$ there exists a polynomial $f(x)\in\Z[x]$ with $\deg f(x)=k-2$ such that there is an $(F_f,2)$-Diophantine set whose size is $k$. More precisely, given a set of integers $S$ containing $k$ distinct elements, we show the existence of a polynomial $f(x)$ with integer coefficients and degree $k-2$ for which $S$ is an $(F_f,2)$-Diophantine set.

Let $f(x)=f_0+f_1 x+\ldots + f_kx^k$. We choose $x_i\in\Z$, $i=0,\ldots, n$, where $x_i\ne x_j$ if $i\ne j$ and $n>k$. We consider the following system of Diophantine equations
\begin{eqnarray}
\label{eq1}
z_{i}^2=f(x_0)f(x_{i}),\;i=1,\ldots,n.
\end{eqnarray}
The variety defined by the latter equations will be called $V_n^k$. In fact, $V_n^k$ lies in the projective space $\PP^{n+k}$ with coordinates $(f_0,f_1,\ldots,f_k,z_1,z_2,\ldots,z_{n})$.
\begin{Lemma}
\label{lem1}
If the system (\ref{eq1}) has a solution $(f_0,f_1,\ldots,f_k,z_1,z_2,\ldots,z_{n})$ with $f_i\in \Z$, $z_i\in\Z$, then $\displaystyle A=\{x_0,x_1,\ldots,x_{n}\}$ is a Diophantine-$(F_f,2)$ set, where $f(x)=f_0+f_1 x+\ldots + f_kx^k$.
\end{Lemma}
\begin{Proof}
If such a solution exists, then it is clear that $f(x_0)f(x_i)$, $i=1,\ldots,n$, is a perfect square. It follows that the product of $f(x_0)f(x_i)$ and $f(x_0)f(x_j)$ is also a perfect square for any $1\le i,j\le n$. Thus, $\displaystyle F(x_i,x_j)=f(x_i)f(x_j)=\left(z_iz_j/f(x_0)\right)^2$, hence a perfect square for any choice of $i$ and $j$.
\end{Proof}

\section{The variety $W_n^k$, $n>k$}

In what follows we are going to write \[\left|\begin{array}{ccccccc}
                 0 & 1 & 2 & 3& \ldots & k& i \\
                 Z_0 & Z_1 & Z_2 & Z_3& \ldots &Z_{k}& Z_i
               \end{array}
\right|,\;i\ge k+1,\] for the $(k+2)\times(k+2)$ determinant
\[\left|\begin{array}{ccccccc}
          1 & 1 & 1 & 1 &\ldots&1&1\\
          x_0 & x_1 & x_2 & x_3 &\ldots&x_{k}&x_i \\
          x_0^2 & x_1^2 & x_2^2 &x_3^2 & \ldots& x_{k}^2&x_i^2 \\
          x_0^3 & x_1^3 & x_2^3 &x_3^3 & \ldots& x_{k}^3&x_i^3 \\
          \vdots &  &   &\vdots   &&& \vdots \\
          x_0^{k} & x_1^{k} & x_2^{k} &x_3^{k} & \ldots& x_{k}^{k}&x_i^{k}\\
          Z_0 & Z_1 & Z_2& Z_3& \ldots&Z_{k}&Z_i
        \end{array}
\right|.\]

Given $x_i\in\Q$, $i=0,1,\ldots,n$, we define the variety $W_{n}^k$ in $\PP^{n}$ by the following system of determinantal equations:
\[\left|\begin{array}{cccccc}
          0 & 1 & 2 & \ldots &k&i\\
          Y_0^2 & Y_1^2 & Y_2^2& \ldots & Y_k^2&Y_i^2
        \end{array}
\right|=0,\;i=k+1,\ldots,n.\]
In fact, the variety $W_{n}^k$ is an intersection of $n-k$ quadrics in $\PP^n$.

\begin{Proposition}
\label{Prop1}
The variety $V_{n}^k$ is $\Q$-birationally equivalent to $W_{n}^k$ for any $n\ge k+1$.
\end{Proposition}
\begin{Proof}
We will give an explicit description of the birational equivalence. We define the map $\phi:V_{n}^k\to W_{n}^k$ by \[(f_0,f_1,\ldots,f_k,z_1,\ldots,z_{n})\mapsto (f(x_0),z_1,\ldots,z_{n}).\] In fact, the point $(f(x_0),z_1,\ldots,z_{n})$ is a point in $W_n^k$ by using the fact that $z_i^2=f(x_0)f(x_i)$ and the linearity of determinants.

We define the map $\psi:W_{n}^k\to V_n^k$ by sending the point $(Y_0,Y_1,\ldots, Y_{n})\in W_{n}^k$ to
{\footnotesize\begin{eqnarray*} \left(\left|\begin{array}{ccccc}
          x_0 & x_1 & x_2 & \ldots & x_k\\
          x_0^2 & x_1^2 & x_2^2 &\ldots & x_k^2\\
          x_0^3 & x_1^3 & x_2^3 &\ldots & x_k^3\\
          \vdots &  &\vdots  & & \vdots\\
          x_0^k & x_1^k & x_2^k & \ldots & x_k^k\\
          Y_0^2 & Y_1^2 & Y_2^2 & \ldots & Y_k^2
        \end{array}\right|, -\left|\begin{array}{ccccc}
          1& 1& 1 & \ldots & 1\\
          x_0^2 & x_1^2 & x_2^2 &\ldots & x_k^2\\
          x_0^3 & x_1^3 & x_2^3 &\ldots & x_k^3\\
          \vdots &  &\vdots  & & \vdots\\
          x_0^k & x_1^k & x_2^k & \ldots & x_k^k\\
          Y_0^2 & Y_1^2 & Y_2^2 & \ldots & Y_k^2
        \end{array}\right|,\left|\begin{array}{ccccc}
          1& 1& 1 & \ldots & 1\\
          x_0 & x_1 & x_2 &\ldots & x_k\\
          x_0^3 & x_1^3 & x_2^3 &\ldots & x_k^3\\
          \vdots &  &\vdots  & & \vdots\\
          x_0^k & x_1^k & x_2^k & \ldots & x_k^k\\
          Y_0^2 & Y_1^2 & Y_2^2 & \ldots & Y_k^2
        \end{array}\right|,\ldots, \right. \\(-1)^k\left.\left|\begin{array}{ccccc}
          1& 1& 1 & \ldots & 1\\
          x_0 & x_1 & x_2 &\ldots & x_k\\
          x_0^2 & x_1^2 & x_2^2 &\ldots & x_k^2\\
          \vdots &  &\vdots  & & \vdots\\
          x_0^{k-1} & x_1^{k-1} & x_2^{k-1} & \ldots & x_k^{k-1}\\
          Y_0^2 & Y_1^2 & Y_2^2 & \ldots & Y_k^2
        \end{array}\right|,DY_0Y_1, DY_0Y_2 ,\ldots,DY_0Y_{n}\right)\end{eqnarray*}}
where \[D=\left|\begin{array}{ccccc}
          1 & 1  & 1 &\ldots&1\\
          x_0 & x_1 & x_2 & \ldots &x_k\\
x_0^2 & x_1^2 & x_2^2 &\ldots & x_k^2\\
          x_0^3 & x_1^3 & x_2^3 &\ldots & x_k^3 \\
          \vdots & & & & \vdots\\
          x_0^k & x_1^k & x_2^k & \ldots & x_k^k
        \end{array}
\right|.\]
This means that the coefficient $(-1)^jf_j$ is obtained by deleting the last column and the $j$-th row in the $(k+2)\times (k+2)$-determinant
\[\left|\begin{array}{cccccc}
          0 & 1 & 2 & \ldots &k&i\\
          Y_0^2 & Y_1^2 & Y_2^2& \ldots & Y_k^2&Y_{i}^2
        \end{array} \right|,\;0\le i\le n, \] where the latter determinant is clearly $0$ if $i\le k$, and it is $0$ by the definition of $W_n^k$ if $i\ge k+1$. Furthermore, the expansion of the determinant is $f(x_i)+(-1)^{k+1}DY_i^2=0$ where $0\le i\le n$. It follows that $f(x_0)f(x_i)=D^2Y_0^2Y_i^2$, hence a point on $V_n^k$.

Checking that $\psi\circ \phi$ and $\phi\circ\psi$ are the identity maps on $V_n^k$ and $W_n^k$ respectively is performed by direct calculations.
\end{Proof}
\begin{Remark}
\label{rem1} We remark that we may assume that the coordinates $Y_i$ on the projective variety $W_n^k$ are integers. Therefore, if $x_i\in\Z$, $i=0,1,\ldots,n$, where $x_i\ne x_j$ when $i\ne j$, then the map $\psi$ will yield a polynomial $f(x)\in\Z[x]$ of degree $k$. In particular, the existence of a rational point on $W_n^k$ together with our choice of integer $x_i$'s imply the existence of an $(F_f,2)$-Diophantine set of size $n+1$, see Lemma \ref{lem1}.
\end{Remark}
\section{The variety $W_{k+1}^k$}

We recall that a variety is {\em rational} if it is $\Q$-birational to $\PP^m$ for some $m\ge 1$. In this section, we fix our choice of $x_i\in\Q$, $i=0,1,\ldots,k+1$, where $x_i\ne x_j$ when $i\ne j$. We then prove that the variety $W_{k+1}^k$ is rational. Hence $V_{k+1}^k$ is also rational, see Proposition \ref{Prop1}.

\begin{Theorem}
\label{thm1}
Let $x_i\in\Q$, $i=0,1,\ldots,k+1$, be such that $x_i\ne x_j$ when $i\ne j$. Then the variety $W_{k+1}^k$ is rational.
\end{Theorem}
\begin{Proof}
Recall that $W_{k+1}^k$ is defined by the determinant \[\left|\begin{array}{cccccc}
          0 & 1 & 2 & \ldots &k&k+1\\
          Y_0^2 & Y_1^2 & Y_2^2& \ldots & Y_k^2&Y_{k+1}^2
        \end{array}
\right|=0.\] Therefore, $W_{k+1}^k$ is a quadric hypersurface in $\PP^{k+1}$ with coordinates $(Y_0,\ldots,Y_{k+1})$. Furthermore, the point $P=(1,1,\ldots,1)$ lies in the set of rational points $W_{k+1}^k(\Q)$. In particular, $W_{k+1}^k$ is birationally equivalent to $\PP^k$. More precisely, we have a rational map $\phi:\PP^k\to W_{k+1}^k$ such that $\phi(Q)$, where $Q=(q_0,q_1,\ldots,q_k,0)$, is the intersection of the line $L$ that joins the points $P$ and $Q$ with $W_{k+1}^k$. For ease of calculations, we assume that the line spanned by $P$ and $Q$ is given by $2\mu Q-\nu P$. Using determinantal properties, one may show that $\phi(Q)$ is described as follows:
\begin{eqnarray*}
(Y_0,Y_1,\ldots,Y_k,Y_{k+1})=(2\mu q_0-\nu,2\mu q_1-\nu,\ldots,2\mu q_k-\nu,-\nu)
\end{eqnarray*}
where
\begin{eqnarray*}
\mu=\left|\begin{array}{cccccc}
                 0 & 1 & 2 &  \ldots & k& k+1 \\
                 q_0 & q_1 & q_2 & \ldots &q_{k}& 0
               \end{array}
\right|,\textrm{ and \hskip20pt}\nu=\left|\begin{array}{cccccc}
                 0 & 1 & 2 &  \ldots & k& k+1 \\
                 q_0^2 & q_1^2 & q_2^2 & \ldots &q_{k}^2& 0
               \end{array}
\right|.
\end{eqnarray*}
One may easily check that the inverse rational map $W_{k+1}^k\to \PP^k$ is described as follows:
\begin{eqnarray*}
(Y_0,Y_1,\ldots,Y_k,Y_{k+1})\mapsto(Y_0-Y_{k+1},Y_1-Y_{k+1},\ldots, Y_k-Y_{k+1}).
\end{eqnarray*}
\end{Proof}
\begin{Corollary}
\label{cor1}
Given $x_i\in\Z$, $i=0,\ldots, k+1$, where $x_i\ne x_j$ if $i\ne j$, there are infinitely many polynomials $f(x)\in\Z[x]$ with $\deg f=k$ such that the set $\displaystyle\{x_i:0\le i\le k+1\}$ is an $(F_f,2)$-Diophantine set where $F_f(x,y)=f(x)f(y)$.
\end{Corollary}
\begin{Proof}
This follows from Proposition \ref{Prop1} since $V_{k+1}^k$ is birational to $W_{k+1}^k$ where the latter is rational, see Theorem \ref{thm1}. Now one concludes using Remark \ref{rem1}.
\end{Proof}

\section{Constructing $f$ of smaller degree}
In Corollary \ref{cor1}, given any subset $S$ of $\Z$ of size $k$ we may find a polynomial $f(x)$, with $\deg f=k-2$, such that $S$ is an $(F_f,2)$-Diophantine set. In this section we investigate the possibility of finding $f'(x)\in\Z[x]$ of smaller degree than $k-2$ such that $S$ is still an $(F_{f'},2)$-Diophantine set.

We prove that given a set $S$ consisting of $k$ integers, there is a polynomial $f(x)\in\Z[x]$ of degree $\displaystyle 2\left\lfloor\frac{k}{3}\right\rfloor$ such that $S$ is an $(F_f,2)$-Diophantine set, and $\lfloor\textrm{ . }\rfloor$ is the greatest integer function.

\begin{Theorem}
\label{thm:Wnrational}
Let $x_i\in\Q$, $i=0,1,\ldots,3k+1$, be such that $x_i\ne x_j$ when $i\ne j$. Then the variety $W_{3k+1}^k$ is rational.
\end{Theorem}
\begin{Proof}
 The variety $W_{3k+1}^k$ is defined by the intersection of the following $k+1$ diagonal quadrics in $\PP^{3k+1}$
 \[\left|\begin{array}{cccccc}
          0 & 1 & 2 & \ldots &2k&i\\
          Y_0^2 & Y_1^2 & Y_2^2& \ldots & Y_{2k}^2&Y_{i}^2
        \end{array}
\right|=0,\textrm{\hskip20pt} i=2k+1,\ldots,3k+1.\]
  The variety $W_{3k+1}^k$ contains the plane $\Pi$ spanned by the following $k+1$ points \begin{eqnarray*}
 T_0=(1,1,1,\ldots,1),\;T_1&=&(x_0,x_1,x_2,\ldots,x_{3k+1}),\;T_2=(x_0^2,x_1^2,x_2^2,\ldots,x_{3k+1}^2),\ldots,\\T_{k}&=&(x_0^k,x_1^k,x_2^k,\ldots,x_{3k+1}^k).
 \end{eqnarray*}
  In order to prove the rationality of $W_{3k+1}^k$, we will project away from $\Pi$. More precisely, we display the following birational map $\phi: \PP^{2k}\to W_{3k+1}^k$ defined by sending $Q=(q_0,\ldots,q_{2k},0,\ldots,0)$ to the intersection of the $(k+1)$-dimensional plane $L$ passing through $Q$ and $\Pi$ with $W_{3k+1}^k$. We may define $L$ by $\mu_0 T_0+\mu_1T_1+\ldots+\mu_{k}T_k+\mu_{k+1} Q$. The map $\phi$ is defined explicitly once we give the values of the $\mu_i$'s.
 An intersection point of $L$ with $W_{3k+1}^k$ is a point which lies on $L$ and the $k+1$ diagonal quadrics defining $W_{3k+1}^k$. Therefore, using the properties of the determinant, one may conclude that an intersection point will be a solution for the following system of equations:
\[A\left(
     \begin{array}{c}
       \mu_0 \\
       \mu_1 \\
       \vdots \\
       \mu_k \\
       \mu_{k+1}
     \end{array}
   \right)=0
\]
where $A=\left(A_{i,j}\right)$ is a $(k+1)\times (k+2)$-full rank matrix whose entries are given as follows:
\begin{eqnarray*}
A_{m-2k,1}&=&2\left|\begin{array}{ccccc}
                 0 & 1 & \ldots & 2k & m \\
                 q_0 & q_1 & \ldots & q_{2k}  &0
               \end{array}
\right|,\; A_{m-2k,2}=2\left|\begin{array}{ccccc}
                 0 & 1 & \ldots & 2k & m \\
                 q_0x_0 & q_1x_1 & \ldots & q_{2k} x_{2k}&0
               \end{array}\right|,\ldots,\\
               A_{m-2k,k+1}&=&2\left|\begin{array}{ccccc}
                 0 & 1  & \ldots& 2k & m \\
                 q_0x_0^k & q_1x_1^k &\ldots& q_{2k}x_{2k}^k &0
               \end{array}
\right|,\; A_{m-2k,k+2}=\left|\begin{array}{ccccc}
                 0 & 1 & \ldots & 2k & m\\
                 q_0^2 & q_1^2 & \ldots &  q_{2k}^2&0
               \end{array}\right|
\end{eqnarray*}
where $m=2k+1,2k+2,\ldots,3k+1$.

Now a straight forward linear algebra exercise shows that
$\mu_j=(-1)^{j}\det(A_j)$, $j=0,1,\ldots,k+1$, where $A_j$ is the matrix $A$ with the $(j+1)$-th column being removed. Therefore, $W_{3k+1}^k$ is birational to $\PP^{2k}$ via the map $\phi:\PP^{2k}\to W_{3k+1}^k$ defined by sending $(q_0:q_1:\ldots:q_{2k})$ to \begin{eqnarray*}(\mu_0&+&\mu_1x_0+\mu_2x_0^2+\ldots+\mu_kx_0^k+\mu_{k+1}q_0,\;\mu_0 +\mu_1x_1+\mu_2x_1^2+\ldots+\mu_k x_1^k+\mu_{k+1}q_1,\;\ldots,\\
\mu_0&+&\mu_1x_{2k}+\mu_2x_{2k}^2+\ldots+ \mu_k x_{2k}^k+\mu_{k+1}q_{2k},\;
\mu_0+\mu_1x_{2k+1}+\mu_2x_{2k+1}^2+\ldots+\mu_k x_{2k+1}^k,\\
\mu_0&+&\mu_1x_{2k+2}+\mu_2x_{2k+2}^2+\ldots+\mu_k x_{2k+2}^k,\;\ldots,\mu_0+\mu_1x_{3k+1}+\mu_2x_{3k+1}^2+\ldots+\mu_k x_{3k+1}^k).
\end{eqnarray*}
\end{Proof}
Using Theorem \ref{thm:Wnrational}, one obtains the following result.
\begin{Corollary}
Let $x_i\in\Z$, $i=0,1,\ldots,3k+1$, be such that $x_i\ne x_j$ when $i\ne j$. Then there are infinitely many polynomials $f(x)\in\Z[x]$ of degree $2k$ such that $\displaystyle\{x_i:0\le i\le 3k+1\}$ is an $(F_f,2)$-Diophantine set.
\end{Corollary}
\begin{Proof}
The proof is a direct consequence of Theorem \ref{thm:Wnrational} together with the fact that $W_{3k+1}^k$ and $V_{3k+1}^k$ are birational, see Proposition \ref{Prop1}. Now we conclude using Remark \ref{rem1}.
\end{Proof}

The variety $V_n^k$ can possibly give rise to hyperelliptic curves whose Jacobian varieties are of high rank.
More precisely, if $(f_0,\ldots,f_k,z_1,\ldots,z_n)\in V_n^k (K)$ where $K=\Q(x_0,\ldots,x_n)$, then the quadratic twist $C^T:f(x_0)y^2=f(x)=f_0+f_1x+\ldots+f_kx^k$ of the hyperelliptic curve $C:y^2=f(x)$ by $f(x_0)$ contains the following rational points in $C^T(K)$
\begin{eqnarray*}
\begin{array}{ccccc}
  P_0=(x_0,1), & P_1=\left(x_1,\frac{z_1}{f(x_0)}\right),& P_2=\left(x_2,\frac{z_2}{f(x_0)}\right), & \ldots, & P_n=\left(x_{n},\frac{z_{n}}{f(x_0)}\right).
\end{array}
\end{eqnarray*}
The point $P_0$ is an obvious point in $C^T$. For the points $P_i$, $i\ge 1$, one observes that $z_i^2=f(x_0)f(x_i)$, hence $\displaystyle f(x_0)\left(z_i/f(x_0)\right)^2=f(x_i)$.

 We recall that the hyperelliptic curve $C^T$ may be embedded in its Jacobian variety $J^T$ via the Albanese map $j: C^T\hookrightarrow J^T$ defined by $\displaystyle j(P)= \left[(P)-(P_0)\right]$, where $[\textrm{ . }]$ denotes the linear equivalence class of the divisor on $C^T$. The point $P_0$ is called the base point of the embedding $j$. Therefore, the points $j(P_i)$, $1\le i\le n$, generate a subgroup $\sum$ in $J^T(K)$ which is generically of positive rank. A natural question would be how large the rank of $\sum$ might be. In fact, infinitely many elliptic curves with rank $7$ over $\Q$ were obtained by using the variety $V_7^4$, see \cite{Yamagishi3}. The authors intend to investigate this question further in future work.

\end{document}